\documentclass[a4paper]{article}

\usepackage[english]{babel}

\usepackage{geometry}

\usepackage{amsmath}
\usepackage{amsthm}
\usepackage{thmtools}
\usepackage{amsfonts,amssymb}
\usepackage{graphicx}
\usepackage{hyperref}
\hypersetup{colorlinks=true, allcolors=blue}
\usepackage{setspace}
\usepackage{siunitx} 
\usepackage{graphicx,pstricks,listings,stackengine}
\usepackage{enumerate} 
\usepackage{authblk}
\usepackage{float} 
\newtheorem{thm}{Theorem}
\newtheorem{prop}{Proposition}
\newtheorem{lem}{Lemma}

\theoremstyle{definition}

\newtheorem{eg}{Example}
\newtheorem{rmk}{Remark}
\newtheorem{algo}{Algorithm}
\title{Convergence, Finiteness and Periodicity of Several New Algorithms of $p$-adic Continued Fractions}
\author{Zhaonan Wang}
\author{Yingpu Deng}
\affil{Key Laboratory of Mathematics Mechanization, NCMIS, Academy of Mathematics and Systems Science, Chinese Academy of Sciences, Beijing 100190, People's Republic of China\authorcr and\authorcr 
	School of Mathematical Sciences, University of Chinese Academy of Sciences,
	Beijing 100049, People’s Republic of China
	\authorcr
	znwang@amss.ac.cn, dengyp@amss.ac.cn}
\date{}
\begin{document}
\maketitle
\bibliographystyle{alpha}
\linespread{1.1}
\begin{abstract}
Classical continued fractions can be introduced in the field of $p$-adic numbers, where $p$-adic continued fractions offer novel perspectives on number representation and approximation. While numerous $p$-adic continued fraction expansion algorithms have been proposed by the researchers, the establishment of several excellent properties, such as the Lagrange's Theorem for classic continued fractions, which indicates that every quadratic irrationals can be expanded periodically, remains elusive. In this paper, we introduce several new algorithms designed for expanding algebraic numbers in $\mathbb{Q}_p$ for a given prime $p$. We give an upper bound of the number of partial quotients for the expansion of rational numbers, and prove that for small primes $p$, our algorithm generates periodic continued fraction expansions for all quadratic irrationals. Experimental data demonstrates that our algorithms exhibit better performance in the periodicity of expansions for quadratic irrationals compared to the existing algorithms. Furthermore, for bigger primes $p$, we propose a potential approach to establish a $p$-adic continued fraction expansion algorithm. As before, the algorithm is designed to expand algebraic numbers in $\mathbb{Q}_p$, while generating periodic expansions for all quadratic irrationals in $\mathbb{Q}_p$.

\end{abstract}
\section{Introduction}
 Continued fractions, a central topic in number theory, serve as an important tool for Diophantine approximation, offering the finest rational approximations achievable. Notably, the algorithm for computing the continued fraction expansion of a real number yields these optimal approximations.A finite real continued fraction precisely corresponds to a rational number. Regarding infinite continued fractions, Euler demonstrated that classical periodic continued fractions exclusively correspond to quadratic irrationals. Furthermore, Lagrange established in \cite{Lagrange} that every quadratic irrational possesses a periodic expansion. In 1940, Mahler \cite{mahler1940} initiated the exploration of continued fractions within the field of $p$-adic numbers $\mathbb{Q}_p$, which motivated authors to find the expansion algorithm that can maintain analogous properties exhibited by real continued fractions. In the twentieth century, Mahler \cite{mahler1934}, Schneider \cite{schneider1970}, Ruban \cite{ruban1970}, Bundschuh \cite{bundschuh1977}, Browkin \cite{browkinI} and several authors had proposed algorithms to generalize the properties of classical continued fractions to the $p$-adic case. However, there is still no suitable algorithm satisfying properties akin to the Lagrange's theorem in classical continued fractions, since the floor function for choosing partial quotients when defining a $p$-adic continued fraction expansion algorithm can be varied. It has been proved that the algorithms defined by Schneider and Ruban cannot guarantee finite expansions for rational numbers \cite{hirsh2011} \cite{laohakosol1985}, and all of the above algorithms cannot provide periodic expansions for quadratic irrationals. In 2001, Browkin introduced four algorithms in \cite{browkinII}, with the first two applicable to all $p$-adic numbers (denoted as Browkin I and Browkin II). Like Browkin I, Browkin II also ensures finite-length expansions for rational numbers, as proven in \cite{barbero2021periodic}. In the case of Browkin I, Bedocchi has characterized the purely periodic expansions and the periodic behavior of certain quadratic irrational numbers in \cite{bedocchi1988nota} and \cite{bedocchi1989remarks}. In \cite{capuano2023periodicity}, Capuano $\textit{et~al.}$ provided additional insights and established more general properties on the periodicity of Browkin I. The other two algorithms proposed in \cite{browkinII} are for quadratic irrationals only, providing better experimental results in representing quadratic irrationals periodically compared to the algorithms in \cite{browkinI}.  \par
In recent years, Nadir Murru and his collaborators have contributed significantly to this topic, discussing the convergence \cite{murru2023convergence} and periodicity properties \cite{murru2023periodicity} of $p$-adic continued fractions. Their works have also generalized the $p$-adic continued fractions in $\mathbb{Q}_p$ to the $\mathfrak{P}$-adic continued fractions in $K $ \cite{murru2021finiteness}, where $\mathfrak{P}$ is a prime ideal in an algebraic number field $K$. Unfortunately, an analogous version of Lagrange's Theorem has yet to be proven in the context of all these algorithms, since the behavior of periodicity of $p$-adic continued fraction expansions relies heavily on the floor functions for the partial quotients, and the diversity in available choices for these floor functions further contributes to the complexity. In \cite{romeo2023continued}, Romeo gave a concise overview of the main results, developments and open problems in the theory of $p$–adic continued fractions. In 2023, Murru $\textit{et~al.}$ introduced a new algorithm \cite{murru2023new},  which modifies and improves one of Browkin's algorithms in \cite{browkinII}, and is considered one of the best at the present time. Although empirical evidence indicates that this novel algorithm indeed yields periodic continued fraction expansions for more purely quadratic irrationals compared to the old ones, a theoretical foundation for this phenomenon is still absent. Since this algorithm chooses floor functions based on the parity of indices (like the algorithms introduced by Browkin in \cite{browkinII}), it leads to an obvious manifestation of this phenomenon in most properties when discussing the periodicity and pure periodicity, including the periodic and pre-periodic lengths. It is our perspective that this phenomenon arises from the characteristics introduced by the algorithm's decisions, rather than being intrinsic to the nature of $p$-adic continued fractions themselves. The fundamental reason behind this lies in the algorithm's methods of choosing partial quotients based on the parity of their indices, which, in turn, leads to disparities in the $p$-adic valuations of partial quotients and complete quotients when the index is odd versus even. For example, when $n$ is odd, it is typically observed that $v_p(\alpha_n)<0$, followed by $v_p(\alpha_{n+1})=0$,  unless the coefficient $a_0$ of index 0 in the $p$-adic expansion of $\alpha_n$ is 0. Similarly, when $n$ is even, $v_p(\alpha_n)$ and $v_p(b_n)$ appear to be 0, consequently leading to $v_p(\alpha_{n+1})<0$. Since periodicity implies $\alpha_{n}=\alpha_{n+k}$ for some $k$, the length of the period will inevitably exhibit certain characteristics of parity.\par
In \cite{barbero2021periodic}, the authors demonstrated that periodic continued fractions with real numbers as partial quotients typically exhibit convergence under the real metric. It's important to note that the current definition methods of $p$-adic continued fraction algorithms still utilize rational numbers as partial quotients. Hence, to obtain periodic $p$-adic continued fractions for quadratic irrationals, there arises a need for an algorithm that integrates the properties of approximating numbers from both $p$-adic and real metric perspectives. Consequently, we contend that, given our current goal to find an algorithm capable of providing periodic expansions for all quadratic irrationals in $\mathbb{Q}_p$, the properties of numbers as real or algebraic numbers must be considered. In other words, merely possessing information about a number's $p$-adic properties (like its $p$-adic expansion) is insufficient to enable us to obtain an algorithm that fulfills the aforementioned objective, this is also why the trace of algebraic numbers is utilized in our algorithms. Furthermore, although the discussion in this paper primarily centers around rational numbers and quadratic irrationals in $\mathbb{Q}_p$, for the sake of comprehensive definition, we present the algorithms applicable to all algebraic numbers in $\mathbb{Q}_p$. \par
In this paper, we recall some concepts, properties, as well as algorithms of $p$-adic continued fractions in Section \ref{sec2}. In Section \ref{sec3}, we firstly introduce two novel algorithms designed to expand algebraic numbers in $\mathbb{Q}_p$. We prove that the algorithms terminate exactly when expanding rational numbers, and give an upper bound for the number of expansions depending only on the denominator of the rational number. Moreover, we prove that, for $p=2,3$, these two algorithms always generate periodic expansions for all quadratic irrationals in $\mathbb{Q}_p$. Inspired by this proof, we propose a potential approach to establish $p$-adic continued fraction expansion algorithms for bigger primes $p$. As before, the algorithms are designed to expand algebraic numbers in $\mathbb{Q}_p$, while generating periodic expansions for all quadratic irrationals in $\mathbb{Q}_p$. We provide a detailed example specifically for $p\leq7$. We present in Section \ref{sec4} some computational results, and finally make a conclusion in Section \ref{sec5}.

\section{Preliminaries}\label{sec2}

In this paper, $p$ is a prime number, and we denote by $v_p(\cdot)$ and $|\cdot|_p$ the $p$-adic valuation and the $p$-adic norm. First we recall some definitions and notations of \cite{browkinI} and \cite{browkinII}.\par
Let $\mathcal{R}\subset \mathbb{Q} $ be a set of representatives modulo $p$ such that $0\in\mathcal{R}$. Then every $\alpha\in\mathbb{Q}_p$ can be written in the form
$$\alpha=\sum_{n=r}^{\infty}a_n p^n,$$
where all $a_n\in\mathcal{R}$, $r=v_p(\alpha)$ and $a_r\neq0$ if $\alpha\neq0$. We can define the functions $s$ and $t$ as follows:
$$s(\alpha)=\sum_{n=r}^0 a_n p^n,\ t(\alpha)=\sum_{n=r}^{-1} a_n p^n.$$
To some extent, these two functions can be regarded as the ``fractional part" of $\alpha$, as both $\alpha-s(\alpha)$ and $\alpha-t(\alpha)$ are $p$-adic integers. In particular, $\mathcal{R}$  can be taken as $\{-\frac{p-1}{2},\cdots,0,\cdots,\frac{p-1}{2}\}$ or $\{0,1,\cdots,p-1\}$, which correspond to the floor function choices of Browkin \cite{browkinI} and Ruban \cite{ruban1970}, respectively. Given the limitation that Ruban's choice of $\mathcal{R}$ fails to ensure finite $p$-adic continued fraction expansions for all rational numbers, we choose Browkin's definition. i.e.,
$a_n\in\left\{0,\pm1,\cdots,\pm\frac{p-1}{2}\right\}$ for all $\alpha\in\mathbb{Q}_p$.
We adopt the notations introduced in \cite{murru2023periodicity}
 and \cite{murru2023new} and define the sets
 $$
J_p=\left\{\frac{a_0}{p^n} \mid n \in \mathbb{N},-\frac{p^{n+1}}{2}<a_0<\frac{p^{n+1}}{2}\right\}=\mathbb{Z}\left[\frac{1}{p}\right] \cap\left(-\frac{p}{2}, \frac{p}{2}\right),
$$
and
$$
K_p=\left\{\frac{a_0}{p^n} \mid n \geq 1,-\frac{p^n}{2}<a_0<\frac{p^n}{2}\right\}=\mathbb{Z}\left[\frac{1}{p}\right] \cap\left(-\frac{1}{2}, \frac{1}{2}\right) .
$$
I.e., $J_p$ is the set of all possible values taken by the function $s$, and $K_p$ is the analogue set for the function $t$.
 \begin{prop}
 Let $|\cdot|$ be the Euclidean norm, then we have:
 \item (1) For all $a,b\in J_p$ and $a\neq b$, $v_p(a-b)\leq0$;
 \item (2) For all $a,b\in K_p$ and $a\neq b$, $v_p(a-b)<0$;
 \item (3) Let $\alpha\in\mathbb{Q}_p$, then $|s(\alpha)|<\frac{p}{2}$ and $|t(\alpha)|<\frac{1}{2}$.
 \end{prop}
 The proof can be found in \cite{bedocchi1988nota}, \cite{mahler1940} and \cite{murru2023periodicity}.\par
Now we recall some definitions and properties related to continued fractions. A continued fraction expansion, whether finite or infinite, can be written in the form 
$$[b_0,b_1,b_2,\cdots] = b_0 + \dfrac{1}{b_1+ \dfrac{1}{b_2+\cdots}},$$
 and call $\{b_n\}$ the partial quotients.
Given an expansion of a $p$-adic number $\alpha = \sum\limits_{n=r}^\infty a_n p^n$ with a partial quotient sequence $\{b_n\}$, one can calculate the numerators $A_n$ and denominators $B_n$ of the convergents $[b_0,\cdots,b_n]$:
$$A_{-1}=1,\ A_0=b_0,\ A_n=b_n A_{n-1}+A_{n-2},\text{ for n}\geq1, $$
$$B_{-1}=0,\ B_0=1,\ B_n=b_n B_{n-1}+B_{n-2},\text{ for n}\geq1, $$
$$\frac{A_n}{B_n}=[b_0,b_1,\cdots,b_n].$$
For all $n\geq0$, we have 
$$A_n B_{n-1}-B_nA_{n-1}=(-1)^{n-1}. $$
Moreover, denote $\alpha_0 = \alpha$, one can define the complete quotients $\{\alpha_n\}$ for $\alpha$ recursively as follows:
$$\alpha_{n+1}= \dfrac{1}{\alpha_n - b_n}.$$
One can see \cite{hensley2006continued} for more properties of continued fractions.\par
 In the work by Browkin \cite{browkinII}, four $p$-adic continued fraction expansion algorithms were introduced. The initial algorithm was defined in his prior publication \cite{browkinI}, while the subsequent three were presented for the first time in the current study. The first two algorithms work for all $p$-adic numbers, and the latter two are specifically designed for the expansion of quadratic irrationals. To enhance clarity, we will exclusively introduce the second and fourth algorithms in this section. This selection is grounded in their direct relevance to the algorithms we aim to propose in the subsequent section. For ease of reference, we shall designate them as Algorithm \ref{Browkin1} and Algorithm \ref{Browkin4}, respectively, in the order of their appearance in this paper.

\begin{algo}\label{Browkin1}
 For a given $\alpha\in\mathbb{Q}_p$, define inductively (finite or infinite) the sequence $\{b_n\}$ as the partial quotients of the $p$-adic continued fraction expansion of $\alpha$:
$$\begin{cases} \alpha_0=\alpha \\ b_n=s\left(\alpha_n\right) & \text { if } n \text { even } \\ b_n=t\left(\alpha_n\right) & \text { if } n \text { odd and } v_p\left(\alpha_n-t\left(\alpha_n\right)\right)=0 \\ b_n=t\left(\alpha_n\right)-\operatorname{sign}\left(t\left(\alpha_n\right)\right) & \text { if } n \text { odd and } v_p\left(\alpha_n-t\left(\alpha_n\right)\right) \neq 0 \\ \alpha_{n+1}=\dfrac{1}{\alpha_n-b_n} . & \end{cases}$$
\end{algo}
\begin{algo}\label{Browkin4}
$\textnormal{(For quadratic irrationals only)}$
Let $D$ be an integer such that $\sqrt{D}\in\mathbb{Q}_p\backslash\mathbb{Q}$. Then every $\alpha\in\mathbb{Q}(\sqrt{D})\backslash\mathbb{Q}$ can be written  in the form $\alpha=\frac{P+\sqrt{D}}{Q}$ with $P,Q\in\mathbb{Q}.$ \par
Define the function $$s_1(\alpha)=s^\prime,$$
where  $s^\prime\in\{s(\alpha),s(\alpha)-p\cdot\operatorname{sign}(s(\alpha))\}$,
such that $$\left|\frac{P}{Q}-s^\prime\right|=\min\left\{\left|\frac{P}{Q}-s(\alpha)\right|,\left|\frac{P}{Q}-(s(\alpha)-p\cdot\operatorname{sign}(s(\alpha)))\right|\right\}.$$
Similarly define 
$$t_1(\alpha)=t^\prime,$$
where $t^\prime\in\{t(\alpha),t(\alpha)-\operatorname{sign}(t(\alpha))\}$, such that
$$\left|\frac{P}{Q}-t^\prime\right|=\min\left\{\left|\frac{P}{Q}-t(\alpha)\right|,\left|\frac{P}{Q}-(t(\alpha)-\operatorname{sign}(t(\alpha)))\right|\right\}.$$
Then $b_n$ is defined recursively as follows:
$$\begin{cases} \alpha_0=\alpha \\ b_n=s_1(\alpha_n) & \text { if } n \text { even } \\ b_n=t_1(\alpha_n) & \text { if } n \text { odd  } \\ \alpha_{n+1}=\dfrac{1}{\alpha_n-b_n} . & \end{cases}$$
\end{algo}
Browkin observed that, for the case $p=5$, the first algorithm applied to $\sqrt{m}$ fails to give periodic 5-adic continued fraction expansion when $m=$19, 26, 29, 31, 39, 41, 44, 46, 51, 56, 59, 61, 66, 71, 79, 84, 86, 89, 91, 96 (20 cases) for $1\leq m \leq 100$,  while the second algorithm gives periodic 5-adic continued fraction expansion  for every $2\leq m \leq 5000$.\par
In 2023, Murru $\textit{et~al.}$ introduced another algorithm that improves Browkin's algorithm  \ref{Browkin1} \cite{murru2023new}. In this algorithm, the use of the sign function is omitted, and they claimed that this omission contributes to a more streamlined formulation, without compromising the inherent advantage of the original algorithm. Moreover, they suggested that this adjustment enhances significantly the periodicity properties. In this way, given $\alpha_0=\alpha$,  for $n\geq 0$, the algorithm works as follows:
\begin{algo}\label{murru}
$$\begin{cases} 
b_n=s(\alpha_n) & \text { if } n \text { even } \\ b_n=t(\alpha_n) & \text { if } n \text { odd  } \\ \alpha_{n+1}=\dfrac{1}{\alpha_n-b_n} . & \end{cases}$$
\end{algo}

Next we present a lemma that gives the $p$-adic valuation of $A_n$ and $B_n$.
\begin{lem}\label{valuation}
    If for all $n\geq 1$, $b_n\neq0,\ v_p(b_n )\leq 0$ and $v_p(b_n b_{n+1} )<0$, then $$\ v_p(B_n)=\sum_{k=1}^n v_p(b_k)\text{ for n}\geq 1,\ v_p(A_n)=\begin{cases}
        \sum\limits_{k=0}^n v_p(b_k) & \text{ for all }n\geq0,\text{ if } b_0\neq0 \\
        \sum\limits_{k=2}^n v_p(b_k) & \text{ for all }n\geq2,\text{ if } b_0=0 
    \end{cases}.$$

\end{lem}
The proof of this lemma is omitted, since similar properties have been demonstrated in several prior works, such as \cite{bedocchi1988nota} and \cite{murru2023convergence}. In fact, among the previously established lemmas, the conditions were more stringent compared to the one presented above.\par
There is a general consensus that a ``good" $p$-adic continued fraction expansion algorithm should adhere to several fundamental criteria: 
\begin{itemize}
\item Convergence in the $p$-adic sense should be assured for the continued fraction for any legal sequence of partial quotients $\{b_n\}$.
\item The algorithm should yield finite expansions for all rational numbers (This is not satisfied by Ruban's approach, see \cite{capuano2019effective}.)

\end{itemize}
Indeed we have
$$
\left|\frac{A_m}{B_m}-\frac{A_n}{B_n}\right|_p=\left|\frac{A_{n+1}}{B_{n+1}}-\frac{A_n}{B_n}\right|_p=\left|\frac{(-1)^n}{B_n B_{n+1}}\right|_p,
$$
if $\lim\limits_{n\rightarrow\infty}v_p(B_n B_{n+1})=-\infty $, then $\left\{\dfrac{A_n}{B_n}\right\}$ is a Cauchy sequence, therefore convergent in $\mathbb{Q}_p$.
The below lemma can be employed to establish the satisfaction of the first criteria of the expansion algorithms.
\begin{lem}\label{lem_murru_conv}
{\rm{\cite{murru2023convergence}}}
Let $b_0,b_1,\cdots\in\mathbb{Z}[\frac{1}{p}]$ be an infinite sequence such that $$v_p(b_n b_{n+1})<0$$
for all $n>0$. Then the continued fraction $[b_0,b_1,\cdots]$ is convergent to a $p$-adic number.
\end{lem}

While Algorithm \ref{murru} has been recognized as the most effective method for achieving periodicity when representing quadratic irrationals compared to the existing algorithms, there remains a noticeable absence of a theoretical basis to explain the performance. Therefore, we continue to seek an algorithm that can produce even more favorable outcomes with respect to the criteria emphasized in this section. 
 
\section{The New Algorithm}\label{sec3}
For any $\alpha\in\mathbb{Q}_p$ with minimal polynomial of degree not exceeding 2, we can write $\alpha=\dfrac{P+\sqrt{D}}{Q}$, where $P\in\mathbb{Z}$, $0\neq Q\in\mathbb{Z}$, and $D\in\mathbb{Z}$ is 0 (when $\alpha\in\mathbb{Q}$) or a non-square integer (not necessarily square-free). It can be observed that, if we denote $P=P_0$, $Q=Q_0$, then all complete quotients $\alpha_n$ can also be written in the form $\dfrac{P_n+\sqrt{D}}{Q_n}$ for $n\geq0$. (Here one should be aware that $P_n,\ Q_n$ may not in $\mathbb{Z}$.) We establish the following lemmas to determine $P_n$ and $Q_n$.
\begin{lem}\label{lem3}
    Given $\alpha=\dfrac{P+\sqrt{D}}{Q}\in\mathbb{Q}_p$ and its $p$-adic continued fraction expansion $[b_0,b_1,\cdots]$ (finite or infinite), assume $v_p(\alpha_n)=-e_n$, $n\geq0$. If $b_n$ is chosen such that $v_p(b_n)=v_p(\alpha_n)=-e_n$ and $v_p(\alpha_n-b_n)\geq0$ for all $n\geq0$, we have the following results:
    \item (1) If $D=0$, then $P_{n+1}=\operatorname{sign}(P_n-b_n Q_n)\cdot q_n$, $Q_{n+1}=\left|\dfrac{P_n-b_n Q_n}{p^{e_n}}\right|$, where 
    $q_n=\begin{cases}
        \frac{Q_n}{p^{e_n}} & e_n\geq0\\
        Q_n & e_n<0
    \end{cases}$.
    Moreover, $P_{n},Q_{n}\in\mathbb{Z}$ for all $n\geq 0$.
    \item (2) {\cite{de1988periodicity}} If $D\neq0$, then $P_{n+1}=b_n Q_n-P_n$, $Q_{n+1}=\dfrac{D-P_{n+1}^2}{Q_n}$.
\end{lem}
\begin{proof}
We prove the first claim by induction. It holds for $n=0$ by definition. Denote $b_n=\begin{cases}
    \dfrac{b_n^\prime}{p^{e_n}} & e_n\geq0\\
    b_n & e_n<0
\end{cases}$ , where $\gcd(b_n^\prime,p)=1$. Assume $\alpha_n=\dfrac{P_n}{Q_n}$ and $\alpha_{n+1}$ is defined, then 
$$\begin{aligned}
\dfrac{1}{\alpha_{n+1}} & =\alpha_n-b_n \\
& =\dfrac{P_n}{Q_n}-\dfrac{b_n^\prime}{p^{e_n}} \\
& =\dfrac{P_n-b_n^\prime q_n}{q_n p^{e_n}}.
\end{aligned}$$
Given that $v_p\left(\dfrac{1}{\alpha_{n+1}}\right)=v_p\left(\dfrac{P_n-b_n^\prime q_n}{q_n p^{e_n}}\right)\geq0$ and $v_p(q_n)\geq0$, it follows that $\dfrac{P_n-b_n^\prime q_n}{p^{e_n}}\in\mathbb{Z}$. Consequently, $$\alpha_{n+1}=\dfrac{q_n}{\dfrac{P_n-b_n^\prime q_n}{p^{e_n}}}.$$
We observe that $b_n^\prime q_n=b_n Q_n$, assuming a positive denominator by default and we have  $$P_{n+1}=\operatorname{sign}(P_n-b_n Q_n)\cdot q_n,\ Q_{n+1}=\left|\frac{P_n-b_n Q_n}{p^{e_n}}\right|.$$
\end{proof}

\begin{lem}\label{lem4}
Assume $\alpha=\dfrac{P+\sqrt{D}}{Q}\in\mathbb{Q}_p$ and $D\neq0$, denote $Q_0=Q=p^{f_0}q_0$, where $\gcd(p,q_0)=1$. Then for all $n\geq 0$, $P_n,Q_n\in\frac{1}{Q}\mathbb{Z}$. If $q_0\ |\ D-P^2$, then $P_n,Q_n\in\mathbb{Z}$ for all $n\geq0$.
\end{lem}
The proof could be found in \cite{miller2007quadratic} page 67, Sec 3.2, Proposition 6 and \cite{capuano2019effective} page 1862-1865, Sec 4, Prop 4.1. Despite using the Schneider and Ruban type continued fraction expansions throughout their proofs respectively, the conclusion hold for our case as well.
\begin{rmk}\label{rmk_selection_of_PQD}
One should be aware that $P,Q,D$ are not uniquely determined. While it is possible to enforce certain constraints like coprimality to ensure the uniqueness, our intention is to facilitate computation and attain more explicit outcomes. Therefore, we do not impose such conditions. Furthermore, we can always require 
$$q_0\ |\ D-P^2,$$
Since this can be achieved by replacing $P,q_0,D$ by $P q_0,q_0^2,D q_0^2$ respectively.
\end{rmk}
We now return to the topic of finding $p$-adic continued fraction expansion algorithms that can generate periodic expansions for quadratic irrationals in $\mathbb{Q}_p$. One should be aware that the $p$-adic continued fraction expansion of a quadratic irrational is periodic if and only if it has finitely many complete quotients, which is equivalent to the finiteness of choices for $P_n$ and $Q_n$. 
\begin{lem}
    \label{lem5}
    Let  $\alpha=\dfrac{P+\sqrt{D}}{Q}$ be a quadratic irrational in $\mathbb{Q}_p$ with a $p$-adic continued fraction expansion $[b_0,b_1,\cdots]$, denote the complete quotients $\alpha_n=\dfrac{P_n+\sqrt{D}}{Q_n}$. Then the following are equivalent:
    \item (1) The expansion is ultimately periodic.
    \item (2) There are only finitely many choices for $\alpha_n$.
    \item (3) $|P_n|$ are upper bounded for all $n$.
    \item (4) $|Q_n|$ are upper bounded for all $n$.
\end{lem}
\begin{proof}
    It suffices to show that (2) is equivalent to (3) and (4), and we only prove the equivalence of (2) and (4) for brevity, as the remaining equivalence can be proved in a similar way. \par
    $(2)\Rightarrow(4)$ is trivial by the above content, and next we deduce (2) from (4). Notice that the finiteness in the selection of $\alpha_n$ is equivalent to the finiteness in the selection of both $P_n$ and $Q_n$, which, in turn, is equivalent to the boundedness of their absolute values due to their discrete nature from Lemma \ref{lem4}. From Lemma \ref{lem3},  $Q_{n+1}=\dfrac{D-P_{n+1}^2}{Q_n}$, and Lemma \ref{lem4} shows that $P_n,Q_n\in\dfrac{1}{Q}\mathbb{Z}$. This implies that all $|P_n|,|Q_n|\geq \dfrac{1}{|Q|}$. Consequently, $|Q_{n+1}|$ is upper bounded, since $|D-P_{n+1}^2|$ is upper bounded by the boundedness of $|P_n|$ and $|Q_n|$ is lower bounded. 
    \end{proof}
Now we present a $p$-adic algorithm working with algebraic numbers that attempts to produce periodic $p$-adic continued fraction expansions for quadratic irrationals, along with finite expansions for rational numbers. The algorithm strategically selects partial quotients that minimize $|P_n|$ when $\alpha$ is a quadratic irrational. For any algebraic number $\alpha\in\mathbb{Q}_p$, we may assume that the minimal polynomial of $\alpha$ over $\mathbb{Q}$ is of degree $n$. Let $\operatorname{Tr}(\alpha)$ represent the trace of $\alpha$ over $\mathbb{Q}$.
\begin{algo}\label{new_alg}
Given an algebraic number $\alpha\in\mathbb{Q}_p$, define the floor functions $\bar{s}$ and $\bar{t}$ as follows:
$$\bar{s}(\alpha)=\operatorname{round}\left(\dfrac{\frac{\operatorname{Tr}(\alpha)}{n}-s(\alpha)}{p}\right)\cdot p+s(\alpha),$$
$$\bar{t}(\alpha)=\operatorname{round}\left(\frac{\operatorname{Tr}(\alpha)}{n}-t(\alpha)\right)+t(\alpha).$$
Here the function $\operatorname{round}(x)$ returns the integer nearest to $x$ under the Euclidean metric.\par
In particular, for $n\leq2$ (which is also the focal point of our discussion in this paper), let $\alpha=\dfrac{P+\sqrt{D}}{Q}\in\mathbb{Q}_p$ as previously defined. The floor functions are then:
$$\bar{s}(\alpha)=\operatorname{round}\left(\dfrac{\frac{P}{Q}-s(\alpha)}{p}\right)\cdot p+s(\alpha),$$
$$\bar{t}(\alpha)=\operatorname{round}\left(\frac{P}{Q}-t(\alpha)\right)+t(\alpha).$$
 In a more detailed context, $\bar{s}(\alpha)=s(\alpha)+m\cdot p$, where the integer $m$ is chosen to minimize $\left|\frac{P}{Q}-(s(\alpha)+m\cdot p) \right|$ for all $m\in\mathbb{Z}$. Similarly, $\bar{t}(\alpha) = t(\alpha) + m $, where $m$ is selected such that $\left|\frac{P}{Q}-(s(\alpha)+m) \right|$ attains its minimum among all integers $m$.\par
Set $\alpha_0=\alpha$, then $\alpha$ can be expanded as follows:
$$\begin{cases}
b_0=\bar{s}(\alpha_0) \\
\alpha_{n+1}=\dfrac{1}{\alpha_n-b_n} & n\geq 0\\
b_n=\bar{s}(\alpha_n) & if\ v_p(\alpha_n)=0 \\
b_n=\bar{t}(\alpha_n) & if\ v_p(\alpha_n)<0
\end{cases}$$
\end{algo}
Recalling again that it is essential for a $p$-adic continued fraction expansion algorithm to satisfy $v_p(\alpha_n)=v_p(b_n)$ at each step, we therefore select the floor functions $\bar{s}$ and $\bar{t}$ such that $\bar{s}(\alpha)-s(\alpha)\in p\mathbb{Z}_p$, and $\bar{t}(\alpha)-t(\alpha)\in\mathbb{Z}_p$. The reason behind taking $\bar{s}(\alpha)-s(\alpha)\in p\mathbb{Z}_p$ rather than in $\mathbb{Z}_p$ is based on the fact that if, for some $n$, $\bar{s}(\alpha_n)-s(\alpha_n)\in\mathbb{Z}_p\backslash p\mathbb{Z}_p$, an undesirable outcome might arise where $v_p(\alpha_k) = 0$ for $k \geq n+1$, causing the $p$-adic expansion to fail to converge.\par
In \cite{miller2007quadratic} the author has also provided an algorithm (called NQCF--New Quadratic Continued Fraction), which follows a Schneider-type approach and also aims to minimize $|P_n|$. Nevertheless, the performance of this algorithm exhibits several drawbacks, primarily due to the limitations inherent in Schneider-type algorithms.\par
In practice, we propose a simplified version of Algorithm \ref{new_alg} to enhance computational efficiency. The difference between the results obtained from the two algorithms is negligible in practical applications, and we will prove that both algorithms exhibit identical desirable properties concerning convergence, finiteness and periodicity.
\begin{algo}\label{neww_alg}
    $$\begin{cases}
        \alpha_0=\alpha \\
        b_n=\bar{s}(\alpha_n) & \text{ if } n \text{ is even}\\
        b_n=\bar{t}(\alpha_n) & \text{ if } n \text{ is odd}\\
        \alpha_{n+1}=\dfrac{1}{\alpha_n-b_n}
    \end{cases}$$
\end{algo}
\begin{rmk}\label{rmk_difference_between_two_algorithms}
    The only distinction between these two algorithms lies in the following criterion: If, for some index $n$, $v_p(\alpha_n)<0$ and $v_p(\alpha_n-\bar{t}(\alpha_n))>0$ (i.e., $b_n = \bar{t}(\alpha_n)$ and the coefficient of $p^0$ in the $p$-adic expansion of $\alpha_n$, denoted as $a_0$, is zero), then Algorithm \ref{neww_alg} will select the floor function of index $n+1$ as $\bar{s}$, while Algorithm \ref{new_alg} will continue to designate  $\bar{t}$ as the floor function. Despite our belief that partial quotients with smaller absolute values might yield better behavior in periodicity, experimental data suggests that the scenario mentioned in this remark occurs remarkably infrequently.
\end{rmk}
\begin{eg}\label{eg_for_two_algs}
    Take $p=5$ and $\alpha = \dfrac{973}{234}$. Then $\alpha_0 = \alpha$ can be represented as a 5-adic expansion 
    $$\dfrac{973}{234} = 2 - 5 - 2\cdot 5^2+\cdots,$$
    hence $s(\alpha_0)=2$, and $$\bar{s}(\alpha_0) = \operatorname{round}\left(\dfrac{\frac{973}{234}-2}{5}\right)\cdot 5 +2 =2.$$
    Then $\alpha_1 = \dfrac{1}{\alpha_0-2} = \dfrac{234}{505}$, which can be expanded as
    $$\dfrac{234}{505}= -\dfrac{1}{5} + 2- 2\cdot 5+\cdots $$
    by both Algorithm \ref{new_alg} and Algorithm \ref{neww_alg} we choose $\bar{t}$ as the floor function, and obtain 
    $$\bar{t}(\alpha_1) = \operatorname{round}\left(\dfrac{234}{505}+\dfrac{1}{5}\right)-\dfrac{1}{5} = 1 - \dfrac{1}{5} = \dfrac{4}{5}. $$
    In the end, we obtain identical expansion result for the two algorithms:
    $$\dfrac{231}{74}=\left[2, \dfrac{4}{5}, -4, \dfrac{7}{5}, -4, \dfrac{3}{5}\right].$$
\end{eg}
\begin{prop}\label{prop_floor_func}
For any $\alpha=\dfrac{P+\sqrt{D}}{Q}\in\mathbb{Q}_p$ with the degree of minimal polynomial not exceeding 2, we have:  
\item (1) $v_p(\alpha-\bar{s}(\alpha))>0$, $v_p(\alpha-\bar{t}(\alpha))\geq 0$.
\item (2) $|Q\cdot \bar{s}(\alpha)-P|<\dfrac{p}{2}|Q|$, $|Q\cdot \bar{t}(\alpha)-P|<\dfrac{1}{2}|Q|$.
\end{prop}
\begin{proof}
    Denote $\alpha = \sum\limits_{n = r}^\infty a_n p^n $, then $$\bar{s}(\alpha) = \sum\limits_{n = r}^0 a_n p^n + k\cdot p,\ \bar{t}(\alpha) = \sum\limits_{n = r}^{-1} a_n p^n + j$$ 
    for some satisfactory integer $k$ and $j$. Therefore, $$\alpha - \bar{s}(\alpha) = \sum\limits_{n=1}^\infty a_n p^n -k\cdot p,\ \alpha - \bar{t}(\alpha) = \sum_{n=0}^\infty a_n p^n - j, $$
    hence the first assertion holds.\par
    To prove the second assertion, it suffices to note that $\left|\bar{s}(\alpha)-\dfrac{P}{Q}\right|<\dfrac{p}{2}$ and $\left| \bar{t}(\alpha)-\dfrac{P}{Q}\right|<\dfrac{1}{2}$ due to the specific selections of $k$ and $j$. 
\end{proof}
In the next subsections we prove that our defined algorithms satisfy the fundamental criteria for a ``good" $p$-adic continued fraction expansion algorithm, i.e., convergence and finiteness for rationals.
\begin{rmk}\label{rmk_common_observation}
    It is worth noting that, both of our proposed algorithms satisfy the condition that, if $b_n = \bar{s}(\alpha_n)$ for some $n$, then $b_{n+1} = \bar{t}(\alpha_{n+1}) $. This observation, in conjunction with Proposition \ref{prop_further_conv}, will play a crucial role in proving those three properties of our algorithms.
\end{rmk}
\subsection{Convergence}
\begin{prop}\label{prop_conv}
The Algorithm \ref{new_alg} and \ref{neww_alg} always produce continued fractions that converge to a $p$-adic number.
\end{prop}
\begin{proof}
It suffices to prove that the two algorithms produce sequences of partial quotients $\{b_i\}$ satisfying $v_p(b_{n}b_{n+1})<0$, then by Lemma \ref{lem_murru_conv} the convergence is established.\par
Observe that the floor functions ensure that $v_p(\alpha_n-b_n)\geq 0$, hence $v_p(b_n)\leq 0$ for all $n\geq 0$. If $v_p(b_n)=0$ for some $n$, then $v_p(\alpha_n)$ must be 0 and $b_n=\bar{s}(\alpha_n)$. Therefore, according to Proposition \ref{prop_floor_func} and Remark \ref{rmk_common_observation}, we deduce that  $b_{n+1} = \bar{t}(\alpha_{n+1
})$ and $v_p(b_{n+1})=-v_p(\alpha_n-b_n)<0$ for both algorithms, and the condition in Lemma \ref{lem_murru_conv} is satisfied.
\end{proof}

\subsection{Finiteness}
\begin{thm}\label{thm_fini}
For $\alpha=\dfrac{P}{Q}\in\mathbb{Q}$, the Algorithm \ref{new_alg} and \ref{neww_alg} stops in $\lceil \ln Q \rceil + 2$ steps. 
\end{thm}
\begin{proof}
Again from Remark \ref{rmk_common_observation}, if $b_n=\bar{s}(\alpha_n)$ for some $n$, then $b_{n+1}=\bar{t}(\alpha_{n+1}) $. 
We claim that for both algorithms, $|P_{n+2}|<\dfrac{|P_{n+1}|}{2}$ and $|Q_{n+2}|<\dfrac{|Q_n|}{4}$ for all $n\geq0$. To maintain clarity in notation, we confine all symbols used below to the context of Algorithm \ref{new_alg}. We consider the expansion of $\alpha$ as $[b_0,b_1,\cdots]$ using Algorithm \ref{new_alg}, and denote the sequence of complete quotients as $\{\alpha_n\}$. We examine two cases $v_p(\alpha_n)=0$ and $v_p(\alpha_n)<0$ separately.\par
If $v_p(\alpha_n)=0$, then $v_p(Q_n)=0$, $b_n=\bar{s}(\alpha_n)$. From Lemma \ref{lem3} we know $$|P_{n+1}|=|Q_n|,\ |Q_{n+1
}|=|\bar{s}(\alpha_n) Q_n-P_n|.$$
Proposition \ref{prop_floor_func} informs that $|Q_{n+1}|<\dfrac{p}{2}|Q_n|$. Since $v_p(\alpha_n-b_n)>0$, we have $b_{n+1}=\bar{t}(\alpha_{n+1})$ by definition. This time $v_p(Q_{n+1})<0$, hence $$|P_{n+2}|=\dfrac{|Q_{n+1}|}{p^{v_p(Q_{n+1})}}<\dfrac{|Q_{n+1}|}{p}<\dfrac{|Q_n|}{2}=\dfrac{|P_{n+1}|}{2},\ |Q_{n+2}|=\dfrac{| P_{n+1}-b_{n+1} Q_{n+1} |}{p^{v_p(Q_{n+1})}}.$$ 
Again Proposition \ref{prop_floor_func} reveals that $$|Q_{n+2}|<\dfrac{|Q_{n+1}|}{2 p^{v_p(Q_{n+1})}}<\dfrac{|Q_{n+1}|}{2p}<\dfrac{|Q_n|}{4}.$$

The proof for $v_p(\alpha_n)<0$ is quite similar. We obtain the result that 
$$|P_{n+1}|=\dfrac{|Q_n|}{p^{v_p(Q_n)}},\ |Q_{n+1}|<\dfrac{|Q_n|}{2 p^{v_p(Q_n)}},$$
hence 
$$|P_{n+2}|\leq|Q_{n+1}|<\dfrac{|P_{n+1}|}{2},\ |Q_{n+2}|<\dfrac{p}{2}|Q_{n+1}|<\dfrac{|Q_n|}{4}.$$
Here we use the symbol $\leq$, as it might happen that $$v_p(\alpha_n)<0 \text{ and } v_p(\alpha_{n+1})=-v_p(\alpha_n-\bar{s}(\alpha_n))<0.$$
Therefore, the claim is proved. By Lemma \ref{lem3},  we know all $P_{n},Q_{n}\in\mathbb{Z}$. Consequently, the sequences of $\{P_n\}$ and $\{Q_n\}$ cannot exhibit an infinite number of terms, since both $|P_n|$ and $|Q_n|$ should be larger than 1, and the number of terms $n$ should satisfy the inequation $$n\leq 2\lceil \log_4 Q\rceil +2 = \lceil \ln Q\rceil +2. $$
\end{proof}
\begin{eg}
    Let us consider the continued fraction expansion of $\dfrac{103}{21}$ and $\dfrac{1309328571134}{235713135135}$ in $\mathbb{Q}_{71}$. Using Algorithm \ref{new_alg} and \ref{neww_alg}, we have
    $$\frac{103}{21}=\left[-12, -\frac{10}{71}, 5\right],$$
   which coincides with the expansion derived from the algorithms presented in \cite{browkinII} and \cite{murru2023new}.
    However, 
     $$\frac{1309328571134}{235713135135}=\left[12683, \frac{19}{71}, 4, -\frac{23}{71}, 40, -\frac{20}{71}, -8, \frac{34}{71}, 13, -\frac{18}{71}, 11, \frac{18}{71}, -26, \frac{3}{71}\right],$$ 
     and the length of expansion steps is 14, while the other two algorithms give the same expansion
     $$\left[-26, -\frac{15}{71}, -1, -\frac{3}{71}, -29, \frac{5}{71}, 31, -\frac{19}{71}, 1, -\frac{22}{71}, 16, \frac{35}{71}, 26, -\frac{8}{71}, 16, \frac{47}{71}, -1\right]$$
     with length 17.  
\end{eg}
\begin{rmk}
   Compared to the rational expansion length bounds established by previous $p$-adic continued fraction expansion algorithms, the upper bound we present exhibits a significant reduction. For instance, in comparison to Algorithm \ref{murru}, this bound has been reduced from $O( |P|+p\cdot|Q| )$ \cite{murru2023new} to $O(\ln |Q|)$. The distinctiveness of our algorithm arises from its consideration of both $p$-adic approximation and, to some extent, Euclidean norm-based approximation.\par
    Though lacking theoretical support, in practical computations, we have observed that the number of expansion steps is significantly smaller than the provided upper bound.
\end{rmk}
The following subsection constitutes the essence of this paper, where we aim to demonstrate the good properties of the newly defined algorithms in terms of periodicity when representing quadratic irrationals in $\mathbb{Q}_p$.

\subsection{Periodicity}
Now we assume $\alpha=\dfrac{P+\sqrt{D}}{Q}\in\mathbb{Q}_p$ and $D\neq0$ is a non-square.
\begin{thm}\label{thm_periodic}
For $p=2,3$, Algorithm \ref{new_alg} and \ref{neww_alg} always expand $\alpha$ as an ultimately periodic continued fraction.
    
\end{thm}
\begin{proof}
    Again, we prove the theorem for only one of the Algorithms \ref{new_alg} and \ref{neww_alg}, since the proof for the other follows a similar pattern. Here we choose the Algorithm \ref{neww_alg}. With the argument in Remark \ref{rmk_selection_of_PQD}, for all $n\geq0$, we may assume $P_n,Q_n\in\mathbb{Z}$ for convenience. Note that this convention only holds within the context of this proof, and might not hold necessarily in practical applications.\par 
    Next we illustrate that for all $n\geq0$, $|Q_n|\leq M:=\max\left\{|Q_0|,|Q_1|,\frac{p^2}{4}|D|+1\right\}$ when $p=2,3$ by induction.\par
    It is evident that the inequation is satisfied by both $Q_0$ and $Q_1$. Let us suppose that It is true for $Q_k$ for all $k\leq 2n-1$, we prove its validity for $k=2n\text{ and }2n+1$. Recall that
    $$P_{n+1}=b_n Q_n-P_n,\ Q_{n+1}=\dfrac{D-P_{n+1}^2}{Q_n}, $$
    hence from Proposition \ref{prop_floor_func} we have 
    $$|P_{2n}|<\dfrac{|Q_{2n-1}|}{2},\ |Q_{2n}|\leq \left|\dfrac{D}{Q_{2n-1}}\right|+\left|\dfrac{P_{2n}^2}{Q_{2n-1}}\right|<\dfrac{|Q_{2n-1}|}{4}+\left|\dfrac{D}{Q_{2n-1}}\right|.$$
   For $|Q_{2n+1}|$, we consider
    $$|P_{2n}|<\dfrac{|Q_{2n-1}|}{2},\ |P_{2n+1}|<\dfrac{p}{2}|Q_{2n}|,$$
    obtaining
    $$|Q_{2n}|\leq \left|\dfrac{D}{Q_{2n-1}}\right|+\left|\dfrac{P_{2n}^2}{Q_{2n-1}}\right|<\dfrac{|Q_{2n-1}|}{4}+\left|\dfrac{D}{Q_{2n-1}}\right|,$$
    and 
    $$|Q_{2n+1} |\leq\left|\dfrac{D}{Q_{2n}}\right|+\left|\dfrac{P_{2n+1}^2}{Q_{2n}} \right|<\dfrac{p^2}{4}|Q_{2n}|+\left|\dfrac{D}{Q_{2n}}\right|.$$
    Therefore,
    $$|Q_{2n+1}|<\dfrac{p^2}{16}|Q_{2n-1}|+\dfrac{p^2}{4}\left|\dfrac{D}{Q_{2n-1}}\right|+\left|\dfrac{D}{Q_{2n}}\right|.$$
    The claim is always valid for $Q_{2n}$, as $|Q_{2n}|<\dfrac{M}{4}+|D|<M$. Let $f(x)=\dfrac{p^2}{4}x+\dfrac{|D|}{x}$, observe that $f(x)$ always takes maximum value on the endpoints of the given internal of $x$. Recall that $|Q_{k}|\geq 1$ for all $k$, thus if $|Q_{2n}|\leq|D|$, then $$|Q_{2n+1}|<\max\left\{\dfrac{p^2}{4}|D|+1,\dfrac{1}{4}+|D|\right\}\leq M.$$
    
    Otherwise when $|Q_{2n}|>|D|$, we observe from $|Q_{2n}| < \dfrac{|Q_{2n-1}|}{4}+\left|\dfrac{D}{Q_{2n-1}}\right|$ that $|Q_{2n-1}| \neq 1$, implying $|Q_{2n-1}| > 2$. Consequently,
    $$|Q_{2n+1}|<\dfrac{p^2}{16}M+\left|\dfrac{p^2D}{12}\right|+1<M,$$
    and by Lemma \ref{lem5} the theorem is proved.
\end{proof}
Although we only establish the proof of periodicity exclusively for continued fraction expansions of quadratic irrationals in $\mathbb{Q}_2$ and $\mathbb{Q}_3$, we conjecture that the periodicity holds true for $p=5$. In Section \ref{sec4} we will illustrate some numerical results in details. \par
 The aforementioned proof offers valuable insights into the sufficient conditions for periodicity when $p$ assumes larger values.
\begin{lem}\label{lem_larger_p}
Let $p\geq5$ be a prime and $\alpha=\dfrac{P+\sqrt{{D}}}{Q}\in\mathbb{Q}_p$ with $D\neq0$, where we assume $Q\ |\ D-P^2$. Let $\left\{\alpha_n=\dfrac{P_n+\sqrt{D}}{Q_n}\right\}$ be an iteration sequence for $\alpha$ and $\{b_n\}$ be the corresponding non-zero sequence of partial quotients such that $\{P_n\}$ and $\{Q_n\}$ satisfy the formulae in Lemma \ref{lem3}. Define $M=\max\left\{|Q|,\dfrac{p^2}{4}|D|+1,\dfrac{4(p^2+1)}{3}\right\}$. If $\{b_n\}$ are recursively defined as follows:
$$\begin{cases}
b_{N k} = \bar{s}(\alpha_{N k}),\\
b_{Nk+1}=\bar{t}(\alpha_{Nk+1}),\\
b_{Nk+2}=\bar{t}(\alpha_{Nk+2}),\\
\cdots\\
b_{Nk+N-1}=\bar{t}(\alpha_{Nk+N-1}),
\end{cases}$$
where $N=\lceil \ln p\rceil$, then for all $n\geq 0$, $|Q_n|\leq \dfrac{p^2}{4} M + 1$.
\end{lem}
\begin{proof}
  From Lemma \ref{lem4}, all $P_n,Q_n\in\mathbb{Z}$. Next we claim that, if $|Q_{k}|\leq M$ and $|Q_{k+1}|>M$ for some $k$, then $|Q_{k+1}|< \dfrac{p^2}{4}M + 1$, and there exists an index $k^\prime\in[k+2,k+N]$ such that $|Q_n|$ decreases strictly for $n\in[k+1,k^\prime]$ and, $|Q_{k^\prime}|<M$, thus proving the lemma.\par
  First we show that, $k$ must be a multiple of $N$, in other words, $b_k=\bar{s}(\alpha_k)$. 
If, on the contrary, $b_k=\bar{t}(\alpha_k)$, then according to Proposition \ref{prop_floor_func}, it follows that $|P_{k+1}|<\left|\dfrac{Q_k}{2}\right|$ and $|Q_{k+1}|<\dfrac{|Q_k|}{4}+\left|\dfrac{D}{Q_k}\right|$, which is consistently less than $M$ since $|Q_{k}|\in[1,M]$. Hence, we conclude that $b_k=\bar{s}(\alpha_k)$. According to Lemma \ref{lem3} and Proposition \ref{prop_floor_func}, we derive $|Q_{k+1}|<\dfrac{p^2}{4}|Q_k|+\left|\dfrac{D}{Q_k}\right|$, also the assumption ensures $|Q_{k+1}|>M$. Therefore, we must have 
$$|Q_k|>D,\ |Q_{k+1}|<\dfrac{p^2}{4} M +1.$$
By the definition of $b_n$, it can be inferred that, for all $n\in[k+2,k+N]$,
$$|P_n|<\dfrac{|Q_{n-1}|}{2},\ |Q_n|<\dfrac{|Q_{n-1}|}{4}+\left|\dfrac{D}{Q_{n-1}}\right|,$$
in particular, 
$$|Q_{k+2}|<\dfrac{|Q_{k+1}|}{4}+\left|\dfrac{D}{Q_{k+1}}\right|<\dfrac{p^2}{16} M+\dfrac{1}{4}+1.$$
If $|Q_{k+2}|<M$, then the claim is satisfied. Otherwise, we proceed our discussion to $Q_{k+3}$ and obtain 
$$|Q_{k+3}|<\dfrac{|Q_{k+2}|}{4}+\left|\dfrac{D}{Q_{k+2}}\right|<|Q_{k+2}|<\dfrac{p^2}{64}+\dfrac{1}{16}+\dfrac{1}{4}+1.$$
Similarly, at each step, we make a decision to either terminate the procedure or proceed with computations. Upon reaching the case $n = k+N$, wherein
$$|Q_{k+N}|<\dfrac{p^2}{4^{N}}M+\sum_{j=0}^{N-1}\dfrac{1}{4^j}<\dfrac{p^2}{4^{N}}M+\dfrac{4}{3},$$
considering the selection of $N$, we have $4^N\geq p^2+1>p^2$. Consequently, the right-hand side of the above inequation is less than $\dfrac{p^2}{p^2+1}M+\dfrac{4}{3}$, which is not greater than $M$ due to the restriction $M\geq \dfrac{4}{3}(p^2+1)$. Thus, the claim is proved.
\end{proof}
 If the condition in the lemma is satisfied, $|Q_n|$ will be bounded, hence this lemma can be perceived as a sufficient condition for achieving periodicity, which opens the door to the potential existence of new algorithms for the expansion of rational and quadratic irrational numbers, such that for larger prime $p$, the algorithms preserve the properties of $p$-adic convergence, finiteness for rationals and periodicity for quadratic irrationals.\par
Specifically, for the cases where $p\leq7$, we illustrate such a satisfactory algorithm.
\begin{algo}\label{modified_alg}
 Given an algebraic number $\alpha\in\mathbb{Q}_p$,
$$\begin{cases}
    \alpha_0=\alpha, \\
    b_{3k}=\bar{s}(\alpha_{3k}), \\
    b_{3k+1}=\bar{t}(\alpha_{3k+1}), \\
    b_{3k+2}=\bar{t}(\alpha_{3k+2}),\\   
    \alpha_{n+1}=\dfrac{1}{\alpha_n-b_n}.
\end{cases}$$
\end{algo}
\begin{eg}
    Again we consider the case in Example \ref{eg_for_two_algs}, i.e., $p=5$ and $\alpha=\dfrac{973}{234}$, but this time we use Algorithm \ref{modified_alg}. The initial two partial quotients remain consistent with those in Example \ref{eg_for_two_algs}, and again we obtain $\alpha_2=-\dfrac{101}{34}.$ However, unlike before, in this instance, we continue to take $\bar{t}$ as the floor function for $\alpha_2$. Since
    $$-\dfrac{101}{34}=1+2\cdot5+\cdots,\ t\left(-\dfrac{101}{34}\right)=0,$$ 
    we have $\bar{t}(\alpha_2)=\operatorname{round}\left(-\dfrac{101}{34}\right)=-3,$
    and $\alpha_3 = \dfrac{1}{\alpha_2+3}=34.$ Consequently, we obtain the 5-adic expansion
    $$\dfrac{973}{234}=\left[2,-\dfrac{4}{5}, -3,34\right].$$
\end{eg}
In the preceding Algorithm \ref{new_alg} and Algorithm \ref{neww_alg}, due to the non-usage of the function $\bar{t}$ within consecutive partial quotient computations, the situation where the function attains a value of 0 does not arise. This is because in both of these algorithms, $t(\alpha_n)$ is always invoked to generate a non-zero integer value, thus eliminating the possibility of resulting in zero through integer addition or subtraction. However, in Algorithm \ref{modified_alg}, the repeated utilization of the $\bar{t}$ function may give rise to a phenomenon (with a great possibility) that $v_p(\alpha_{3k+1}-\bar{t}(\alpha_{3k+1}))=0$, then
$$t(\alpha_{3k+2})=0,\ \bar{t}(\alpha_{3k+2})=\operatorname{round}\left(\dfrac{\operatorname{Tr}(\alpha_{3k+2})}{\deg(\alpha)}\right), $$
when $\left|\dfrac{\operatorname{Tr}(\alpha_{3k+2})}{\deg(\alpha)}\right|<\dfrac{1}{2}$, $\bar{t}(\alpha_{3k+2})=0$, then $\alpha_{3k+3}=\alpha_{3k+1}-\bar{t}(\alpha_{3k+1})$, $v_p(\alpha_{3k+3})=0$ and $b_{3k+3}=\bar{s}(\alpha_{3k+3})\neq0 $. \par
Since it has become customary in continued fraction expansions to avoid the presence of 0 as a partial quotient, when writing the expansion sequence, we can conveniently use the formula
$$[\cdots,b_{n-1},b_n,0,b_{n+1},b_{n+2},\cdots ]=[\cdots,b_{n-1},b_{n}+b_{n+1},b_{n+2},\cdots]$$
to eliminate the occurrence of 0. Nevertheless, due to the variable floor function selection based on the index's modulus 3 in our algorithm, if the occurrence of 0 does indeed arise during the expansion process, we retain it temporarily and apply the aforementioned formula collectively after the whole expansion is finished. For the sake of convenience, in the following contents, we will choose to consider either the original sequence or the sequence obtained by removing the zeros, depending on the specific circumstances. \par
Following a similar process as before, we systematically proceed to establish the properties of convergence, finiteness, and periodicity for Algorithm \ref{modified_alg}. It is noteworthy that the convergence and finiteness hold for all primes $p$.\par
The subsequent lemma is a synthesis of Theorem 9 and Corollary 10 from \cite{murru2023convergence}.
\begin{lem}\label{lem_further_conv}
    Let $b_0, b_1, \ldots \in \mathbb{Q}_p$ such that, for all $n \in \mathbb{N}$ :
$$
 \begin{cases}
v_p (b_{3 n+1} )<0, \\
v_p (b_{3 n+2} )=0, \\
v_p (b_{3 n+3} )=0 .
\end{cases}
$$
If $v_p (b_{3 n+3} b_{3 n+2}+1 )=0$ for all $n \in \mathbb{N}$, then $v_p\left(B_{3 n-2}\right)=v_p\left(B_{3 n-1}\right)=v_p\left(B_{3 n}\right)>v_p\left(B_{3 n+1}\right)$, hence $\lim\limits_{n\rightarrow\infty}v_p(B_n)=-\infty$, and the continued fraction $[b_0,b_1,\cdots]$ is convergent to a $p$–adic number.
\end{lem}

\begin{prop}\label{prop_further_conv}
Consider the sequence of partial quotients after eliminating the element 0. The Algorithm \ref{modified_alg} satisfies that, for all values of $n$ satisfying $v_p(b_n)<0$, $v_p(b_{n+1})=0$ and $v_p(b_{n+2})=0$, we have $v_p (b_{n+1} b_{n+2}+1 )=0$ .
\end{prop}
\begin{proof}
   As the convergence is not affected by the beginning finite terms of the partial quotients, we may assume that $n\neq 0$. Therefore, if $v_p(b_n)<0$, $v_p(b_{n+1})=0$ and $v_p(b_{n+2})=0$ for some $n$, it must happen that
   $$b_n=\bar{t}(\alpha_n),\ b_{n+1}=\bar{t}(\alpha_{n+1})=\operatorname{round}\left(\dfrac{\operatorname{Tr}(\alpha_{n+1})}{\deg(\alpha)}\right)\neq0,\text{ and }b_{n+2}=\bar{s}(\alpha_{n+2}). $$
  Assume $\alpha_{n+1}=a_0+a_1p+\cdots$ is in $\mathbb{Z}_p$ with $a_0\neq0$, since $\alpha_{n+2}$ is also in $\mathbb{Z}_p\backslash p\mathbb{Z}_p$, it follows that $$\alpha_{n+2}\equiv (a_0-b_{n+1})^{-1} \bmod p,\ \bar{s}(\alpha_{n+2})\equiv (a_0-b_{n+1})^{-1} \bmod p, $$ 
hence $b_{n+1}b_{n+2}\equiv b_{n+1}(a_0-b_{n+1})^{-1}\bmod p $. Since
$$b_{n+1}(a_0-b_{n+1})^{-1}=(a_0-b_{n+1})^{-1}(b_{n+1}-a_0+a_0)\equiv -1+a_0(a_0-b_{n+1})^{-1} \bmod p,$$
if $p\ |\ b_{n+1}b_{n+2}+1$, then 
$$p\ |\  a_0(a_0-b_{n+1})^{-1},$$ leading to a contradiction. 
\end{proof}
\begin{thm}\label{thm_overall_p_5_7}
For all prime numbers $p$, Algorithm \ref{modified_alg} always produces continued fractions that converge to a $p$-adic number; When expanding rational numbers, the algorithm will stop in $\lceil \ln Q \rceil + 2$ steps. Specifically, for $p\leq7$, the algorithm always expands a quadratic irrational in $\mathbb{Q}_p$ as an ultimately periodic continued fraction.
\end{thm} 
\begin{proof}
Considering the sequences of partial quotients derived from the application of Algorithm \ref{modified_alg}, we observe that, for any natural number $n$ such that $v_p(b_n)=0$, the following must hold:
$$v_p(b_{n+1})<0,\text{ or }v_p(b_{n+1})=0\text{ and }v_{p}(b_{n+2})<0. $$
From the argument presented in \cite{murru2023convergence}, it becomes apparent that, in order to establish convergence for our algorithm, one should demonstrate that $v_p(B_n)\rightarrow-\infty$, and it suffices to prove that the algorithm satisfies the condition of Lemma \ref{lem_further_conv} under the second case, which is already proved in Proposition \ref{prop_further_conv}. Thus, Algorithm \ref{modified_alg} ensures that each expansion is convergent in a $p$-adic way.\par
The proof of the finiteness property for expanding rational numbers can be established analogously to that of Theorem \ref{thm_fini}. In fact, consider the original sequence of partial quotients obtained from Algorithm \ref{modified_alg} (without eliminating 0), and observe that 
$$|Q_{3k}|<\dfrac{p}{2}|Q_{3k-1}|,\ |Q_{3k+1}|<\dfrac{|Q_{3k}|}{2p},\text{ and }|Q_{3k+2}|<\dfrac{|Q_{3k+1}|}{2}, $$
hence $|Q_{3k+2}|<\dfrac{|Q_{3k-1}|}{8}$, and Algorithm 
 \ref{modified_alg} will again stop in $O(\ln p)$ steps.\par
 Finally we show that for $p=2,3,5,7$, Algorithm \ref{modified_alg} will generate periodic expansions for all quadratic irrationals in $\mathbb{Q}_p$. When $p\leq7$, $N=\lceil \ln p\rceil-1\leq2$, by the definition of Algorithm \ref{modified_alg}, the original sequence of partial quotients $\{b_n\}$ satisfies the condition in Lemma \ref{lem_larger_p}. Therefore, 
 all $|Q_n|$ are bounded, and the expansion is ultimately periodic.
\end{proof}
It can be seen that, finding an algorithm for generating periodic continued fraction expansions for larger $p$ is not a hard task, since one could simply choose more partial quotients $b_n$ as $\bar{t}(\alpha_n)$. However, the challenge lies in ensuring that such an algorithm still satisfies the requirement of $p$-adic convergence. In Section 6 of \cite{murru2023convergence}, the author proved that, if there exists an $r\in\mathbb{N}_+$ such that for all $n\in\mathbb{N}$, 
$$\begin{cases}
    v_p(b_{rn+1})<0 \\
    v_p(b_{rn+i})=0,\ \forall i\in\{2,\cdots,r\},
\end{cases}$$
suppose that for all $n\in\mathbb{N}$,
$$v_p(U_{rn+2}^{(i)})=0\text{ for all }i\in\{2,\cdots,r-1\}, \text{ where }r\geq3,$$
$$v_p(U_{rn+3}^{(i)})=0\text{ for all }i\in\{2,\cdots,r-2\}, \text{ where }r\geq4.$$
Then the continued fraction $[b_0,b_1,\cdots]$ is convergent to a $p$-adic number. Here the family of sequences $U_m^{(n)}$ is defined as follows:
$$
U_m^{(0)}=1, \quad U_m^{(1)}=b_m, \quad U_m^{(n+1)}=b_{m+n} U_m^{(n)}+U_m^{(n-1)}.
$$
As mentioned earlier, the most straightforward approach to define new algorithms is to increase the proportion of $\bar{t}(\alpha_n)$ among all the partial quotients. Therefore, we may consider the following definition in a natural way:
$$\begin{cases}
    \alpha_0=\alpha, \\
    b_{rk}=\bar{s}(\alpha_{rk}), \\
    b_{rk+1}=\bar{t}(\alpha_{rk+1}), \\
    \cdots\\    
    b_{rk+r-1}=\bar{t}(\alpha_{rk+r-1}),\\   
    \alpha_{n+1}=\dfrac{1}{\alpha_n-b_n}.
\end{cases}$$
If for a certain $k$, the partial quotients $b_{rk},\cdots,b_{rk+r-1} $ generated by this algorithm are not zero, then
$$\begin{cases}
    v_p(b_{rk+1})<0 \\
    v_p(b_{rk+i})=0,\ \forall i\in\{2,\cdots,r\}.
\end{cases}$$
It remains to verify whether this algorithm meets the above condition, though.
\section{Computational Results}\label{sec4}
In this section, we collect some experimental data about the performance of periodicity of our algorithms \ref{neww_alg} and \ref{modified_alg} ( results from Algorithm \ref{new_alg} are nearly distinguishable from those generated by Algorithm \ref{neww_alg} ). Also, we revisit some numerical results given in \cite{murru2023new} for the Algorithm \ref{Browkin1} and \ref{murru} to make a comparison. The data indicates that our algorithms performs better in terms of periodicity.\par
Following the computation in \cite{murru2023new}, here we consider the periodicity of the expansions of $p$-adic numbers that are square roots of integers. 
\begin{eg}
    We consider $p=5$. First take $D=19$. As calculated in \cite{murru2023new}, in $\mathbb{Q}_5 $, no period has been detected using Algorithm \ref{Browkin1} up to 1000 steps, while $\sqrt{19}$ has periodic continued fraction using Algorithm \ref{murru}, with expansion 
    $$\sqrt{19}=\left[2,\overline{-\frac{2}{5}, 2, \frac{1}{5}, -2, -\frac{2}{5}, -\frac{12}{5}, \frac{2}{5}, -2, \frac{8}{25}, 2, \frac{1}{5}, -1, -\frac{2}{5},-\frac{8}{5}, \frac{2}{5}, -2, \frac{12}{25}, 2, \frac{2}{5}, -1}\right],$$
    which has pre-period length 1 and period length 20. However, with our Algorithm \ref{neww_alg}, 
    $$\sqrt{19}=\left[2, \overline{\frac{3}{5}, -2, \frac{1}{5}, -3, \frac{2}{5}, -1}\right],$$
    the pre-period length is 1 and period length is 6. Moreover, employing Algorithm \ref{modified_alg} we obtain another periodic expansion 
    $$\sqrt{19}=\left[2,\frac{3}{5},-4,-2,,\cdots,-\frac{6}{5},\overline{1,-2,\frac{1}{25},\frac{14}{5}}\right],$$
    with the pre-period length 11 and period length 4.\par
   For $D=235032571341$, Algorithm \ref{murru} fails to generate a periodic expansion when observing the first 10000 steps. However, $\sqrt{D}$ can be expanded periodically by Algorithm \ref{neww_alg}, the pre-period length is 71 and period length is 810. Alternatively, under Algorithm \ref{modified_alg} the pre-period length is 213 and period length is 1986.\par
   When both the  Algorithm \ref{murru} and our new algorithm produce periodic expansions for a certain $D$, it appears that the expansions generated by the new algorithm tend to have both shorter pre-period and period length. This phenomenon is not coincidental; in our observation, this tendency is almost always observed for various $p$.
\end{eg}
In the Table 1, we show the computational results of the periodicity performance of the algorithms \ref{Browkin1}, \ref{murru} and our algorithms \ref{neww_alg}, \ref{modified_alg}. For clarity, we will refer to these algorithms as ``Browkin", ``Murru", ``Our Algorithm \ref{neww_alg}" and ``Our Algorithm \ref{modified_alg}" in the table. We consider the odd primes $p$ less than 100, and the expansion of $\sqrt{D}$ where $1\leq D\leq1000$, where $D$ is not a square and $\sqrt{D}\in\mathbb{Q}_p$. The table collects the data including the numbers continued fraction expansions of square roots which are periodic within 1000 steps using Murru and Browkin's algorithms (which have been calculated in \cite{murru2023new}), as well as the Algorithm \ref{modified_alg}. As for our Algorithm \ref{neww_alg}, the numbers are chosen such that the expansions of square roots show periodicity within 5000 steps for $p\leq31$, and within 1000 steps for $37\leq p\leq97$. In the last column, we list the total number of the calculated $\sqrt{D}$ for each $p$. All these computations have been carried out using Sagemath.
\begin{table}[htbp]\centering\label{table1}
\begin{tabular}{|c|c|c|c|c|c|}
\hline $p$ & \text { Murru } & \text { Browkin } & \text{Our Algorithm \ref{neww_alg}} & \text{Our Algorithm \ref{modified_alg}} & \text { Total } \\
\hline 3 & 42 & 68 & 345 & 345 & 345 \\
\hline 5 & 81 & 67 & 386 & 386 & 386 \\
\hline 7 & 88 & 68 & 407 & 406 & 407 \\
\hline 11 & 99 & 80 & 428 & 426 & 428 \\
\hline 13 & 106 & 89 & 434 & 433 & 434 \\
\hline 17 & 118 & 109 & 441 & 439 & 441 \\
\hline 19 & 108 & 97 & 445 & 443 & 445 \\
\hline 23 & 113 & 102 & 450 & 448 & 450 \\
\hline 29 & 123 & 111 & 449 & 447 & 453 \\
\hline 31 & 133 & 118 & 443 & 454 & 456 \\
\hline 37 & 125 & 121 & 309 & 454 & 456 \\
\hline 41 & 121 & 117 & 319 & 454 & 457 \\
\hline 43 & 122 & 117 & 297 & 455 & 458 \\
\hline 47 & 117 & 110 & 279 & 453 & 461 \\
\hline 53 & 120 & 118 & 261 & 456 & 460 \\
\hline 59 & 126 & 121 & 273 & 458 & 461 \\
\hline 61 & 133 & 124 & 268 & 460 & 462 \\
\hline 67 & 124 & 121 & 264 & 460 & 462 \\
\hline 71 & 119 & 119 & 261 & 461 & 465 \\
\hline 73 & 128 & 125 & 258 & 459 & 462 \\
\hline 79 & 123 & 120 & 251 & 464 & 468 \\
\hline 83 & 122 & 122 & 242 & 458 & 464 \\
\hline 89 & 131 & 127 & 241 & 463 & 466 \\
\hline 97 & 138 & 135 & 239 & 461 & 464 \\
\hline
\end{tabular}
\caption{Periodicity}
\end{table}
It can be seen that the Algorithm \ref{Browkin1} and \ref{murru} exhibit comparable behavior when applied to large primes, producing periodic expansions for approximately only a quarter of the values in set the set of $D$, while our algorithms show obvious advantage in demonstrating periodicity.\par
For the results generated from Algorithm \ref{neww_alg}, we have substantiated through computation that, many instances of $D$ which do not exhibit periodicity within the first 3000 steps can indeed reveal periodic behavior when observed over a greater number of steps.
For example, we consider the case $p=29$, and $D=354$, employing the Algorithm \ref{Browkin1} and \ref{murru}, we did not observe any periodicity within 10,000 steps. However, with Algorithm \ref{neww_alg}, we can obtain a periodic expansion of $\sqrt{D}$ with pre-period length 6 and period length 9157.\par
Furthermore, it is noteworthy that the Algorithm \ref{modified_alg} even demonstrates superior performance in terms of periodicity compared to the Algorithm \ref{neww_alg}. Most $\sqrt{D}$ exhibit periodicity within the initial 1000 steps of the expansion obtained from applying Algorithm \ref{modified_alg}, while it cannot be achieved by the Algorithm \ref{neww_alg}. As previously shown in the table's preamble, for $p$ greater than 31, we examine the periodicity of the Algorithm \ref{neww_alg} within the first 1000 steps of the expansion, while for $p\leq31$, we consider the initial 5000 steps. Consequently, as $p$ transitions from 31 to 37, a noticeable decrease in the numerical observations becomes evident from the table, while the performance of the Algorithm \ref{modified_alg} remains stable. 
\begin{eg}
Take $p=67$, $D=135$, applying Algorithm \ref{neww_alg} we obtain
$$\sqrt{215}=\left[9,\overline{\frac{9}{67},18}\right],$$
the pre-period length is 1 and the period length is 2, while 
$$\sqrt{215}=\left[9,\overline{\frac{9}{67},9,15,-\frac{28}{67},-2,-13,-\frac{50}{67},
3,4,-\frac{9}{67},-9,-15,-\frac{28}{67},2,13,\frac{50}{67},-3,-4}\right]
$$
using Algorithm \ref{modified_alg}, with pre-period length 1 and period length 18, which has larger period length. \par
However, take $D=4518$, use Algorithm \ref{neww_alg} and \ref{modified_alg}, we obtain
$$\sqrt{4518}=\left[30,\overline{\frac{16}{67},6,-\frac{9}{67},58,\frac{3}{67},-14,\frac{97}{67},32,\frac{3}{67},-18,\frac{29}{67},6,-\frac{7}{67},194 }\right]$$
and 
$$\sqrt{4518}=\left[30,\overline{\frac{16}{67},-4,-20,\frac{32}{67},-2,-40}\right],$$
this time the expansion of Algorithm \ref{modified_alg} has slightly shorter period length.
\end{eg}
At present, we do not yet have a clear understanding of how these two algorithms differ in terms of the period length when both of them can expand a quadratic irrational periodically. However, we find that the Algorithm \ref{modified_alg} performs much better in representing ``large" quadratic irrationals as a periodic continued fraction expansion. For example, take $p=31$, $D=252502352$, applying Algorithm \ref{modified_alg} we can expand $\sqrt{D}$ as a periodic continued fraction with pre-period length 292 and period length 528, while no periodicity can be observed within the first 10000 steps in the expansion of Algorithm \ref{neww_alg}.
\section{Conclusions}\label{sec5}
By summarizing the existing $p$-adic continued fraction expansion algorithms, we introduce the Algorithm \ref{new_alg}, and simplify its steps into the Algorithm \ref{neww_alg} for computational convenience.  We establish the properties of $p$-adic convergence and finiteness of rational number expansions for these new algorithms. Furthermore, we demonstrate its capacity to produce periodic expansions for quadratic irrationals when $p<5$. During this proof process, we gain new insights and propose the Algorithm \ref{modified_alg}. We prove that this algorithm maintains the favorable properties of previous algorithms and can extend the periodicity for expanding quadratic irrationals to $p<8$. Despite the absence of theoretical proofs, experimental data indicate that our newly defined algorithms exhibit nice performance even with larger prime numbers. Remarkably, as evident from the table in Section \ref{sec4}, when examining the periodicity behavior on quadratic irrational expansions, Algorithm \ref{modified_alg} demonstrates the most favorable performance.\par
However, some issues remain unresolved. As the algorithms proposed in this paper are solely applicable to algebraic numbers, it prompts further consideration as to whether there exists an algorithm that can operate on all $p$-adic numbers, possessing the three favorable properties described in this paper. Further investigations are also required to ascertain these properties of our algorithms with larger prime numbers. 
\section*{Acknowledgements}
This work is supported by National Key Research and Development Program of China(No. 2020YFA0712300) and National Natural Science Foundation of China(No. 12271517).

\bibliography{ref}
\end{document}